\documentclass{article}
\usepackage{amsmath}
\usepackage{amsthm}
\usepackage{amsfonts}
\usepackage{amscd}
\let\oldlabel=\label
\def\prellabel{\marginparsep=1em\marginparwidth=44pt
    \def\label##1{\oldlabel{##1}\ifmmode\else\ifinner\else
         \marginpar{{\footnotesize\ \\ \tt
                    ##1}}\fi\fi}}

\let\epsilon=\varepsilon
\let\phi=\varphi
\let\kappa=\varkappa
\def\email#1{{\tt #1}}
\def\SetSize{\fontsize{12}{14.4}\selectfont}

\def\Messina{
\textwidth=30pc
\textheight=48pc
\oddsidemargin=1.5cm
\evensidemargin=1.5cm
\topmargin=2.2cm
\parindent=1pc
\pagestyle{empty}
\def\SetSize{\fontsize{10}{12}\selectfont}
}

\Messina

\makeatletter

\renewcommand\section{\@startsection {section}{1}{\z@}%
                                   {-3.5ex \@plus -1ex \@minus -.2ex}%
                                   {2.3ex \@plus.2ex}%
                                   {} }
\renewenvironment{thebibliography}[1]
     {\vspace{3.5ex \@plus 1ex \@minus .2ex}
      \noindent REFERENCES\par\vspace{2.3ex \@plus.2ex}%
      \list{\@biblabel{\@arabic\c@enumiv}}%
           {\settowidth\labelwidth{\@biblabel{#1}}%
            \leftmargin=\labelwidth
            \@openbib@code
            \usecounter{enumiv}%
            \let\p@enumiv\@empty
            \renewcommand\theenumiv{\@arabic\c@enumiv}}%
      \sloppy
      \clubpenalty4000
      \@clubpenalty \clubpenalty
      \widowpenalty4000%
      \sfcode`\.\@m}
     {\def\@noitemerr
       {\@latex@warning{Empty `thebibliography' environment}}%
      \endlist}
\renewcommand{\@biblabel}[1]{#1.}
\makeatother
\newtheoremstyle{theorem}
  {}
  {}
  {\itshape}
  {}
  {}
  {.}
  {.5em}
  {}

\newtheoremstyle{definition}
  {}
  {}
  {}
  {}
  {\MakeUppercase}
  {.}
  {.5em}
  {}

\theoremstyle{theorem}
\newtheorem{theorem}{THEOREM}
\newtheorem{lemma}[theorem]{LEMMA}
\newtheorem{corollary}[theorem]{COROLLARY}
\newtheorem{proposition}[theorem]{PROPOSITION}

\theoremstyle{definition}
\newtheorem{remark}[theorem]{REMARK}

\let\Cal=\mathcal
\let\frak=\mathfrak
\let\Bbb=\mathbb
\let\phi=\varphi

\def\CC{{\NZQ C}}
\def\pp{{\frak p}}
\def\Cok{\operatorname{Cok}}
\def\grade{\operatorname{grade}}
\def\Hom{\operatorname{Hom}}
\def\Ext{\operatorname{Ext}}
\def\rank{\operatorname{rank}}

\def\Ker{\operatorname{Ker}}
\def\Im{\operatorname{Im}}
\def\I{\operatorname{I}}
\def\SS{\operatorname{S}}
\def\PP{{\Bbb P}}

\def\KK{{\Cal K}}

\let\oldbigwedge\bigwedge
\def\BIGwedge{{\textstyle\oldbigwedge}}
\def\medwedge{{\scriptstyle\oldbigwedge}}
\def\bigwedge{\mathchoice{\BIGwedge}{\BIGwedge}{\medwedge}{}}

\let\hat=\widehat
\def\BB{{\operatorname{\Cal B}}}
\def\KK{{\operatorname{\Cal K}}}
\def\RR{{\operatorname{\Cal R}}}
\def\CC{{\operatorname{\Cal C}}}
\let\tensor=\otimes
\let\iso=\cong
\def\projdim{\operatorname{projdim}}
\let\epsilon=\varepsilon
\let\tilde=\widetilde

\begin{document}
\thispagestyle{empty} \vspace*{1.5in} {\fontsize{14}{16.8}\selectfont \noindent
\uppercase{The Koszul complex in projective\\ dimension one}\par}

\SetSize \vspace{2\baselineskip}

\noindent\uppercase{Winfried Bruns}, Universit\"at Osnabr\"uck, FB
Mathematik/Informatik, 49069 Osna\-br\"uck, Germany,
\email{Winfried.Bruns@mathematik.uni-osnabrueck.de}\\[1\baselineskip]
\noindent\uppercase{Udo Vetter}, Universit\"at Oldenburg, FB Mathematik, 26111
Oldenburg, Germany, \email{vetter@\allowbreak mathematik.uni-oldenburg.de}

\vspace{4\baselineskip plus 1 \baselineskip minus 1\baselineskip}

Let $R$ be a noetherian ring and $M$ a finite $R$-module. With a linear form
$\chi$ on $M$ one associates the Koszul complex $K(\chi)$. If $M$ is a free
module, then the homology of $K(\chi)$ is well-understood, and in particular it
is grade sensitive with respect to $\Im\chi$.

In this note we investigate the case of a module $M$ of projective dimension
$1$ (more precisely, $M$ has a free resolution of length $1$) for which the
first non-vanishing Fitting ideal $\I_M$ has the maximally possible grade
$r+1$, $r=\rank M$. Then $h=\grade \Im\chi\le r+1$ for all linear forms $\chi$
on $M$, and it turns out that $H_{r-i}(K(\chi))=0$ for all even $i<h$ and
$H_{r-i}(K(\chi))\iso \SS^{(i-1)/2}(C)$ for all odd $i<h$ where $\SS$ denotes
symmetric power and $C=\Ext_R^1(M,R)$, in other words, $C=\Cok\psi^*$ for a
presentation
$$
0\to F\stackrel{\psi}{\to} G \to M\to 0.
$$
Moreover, if $h\le r$, then $H_{r-h}(K(\chi))$ is neither $0$ nor isomorphic to
a symmetric power of $C$, so that it is justified to say that $K(\chi)$ is
grade sensitive for the modules $M$ under consideration.

We furthermore show that the maximally possible value $\grade \Im\chi=r+1$ can
only occur in two extreme cases: (i) $r=1$ or (ii) $\rank F=1$ and $r$ is odd.

The note was motivated by a result of Migliore, Nagel, and Peterson (see [MNP],
Proposition 5.1). They implicitly prove the result on $K(\chi)$ for Gorenstein
rings $R$, using local cohomology. Our method allows more general assumptions.
(Even the assumption that $R$ is noetherian is superfluous if one uses the
correct notion of grade.) It is based on results in [BV1] and has a predecessor
in [HM]. The case in which $\rank F=1$ has been treated in [BV3].

The situation in which the Fitting ideal $\I_M$ of $M$ has only grade $r$ is
also of interest. For example, it occurs for the K\"ahler differentials of
complete intersections with isolated singularities. While our method also
yields results in this case, we have restricted ourselves to the case of grade
$r+1$ for the sake of clarity.

The detailed account of the linear algebra of $M$ and its exterior powers has
been given in [BV2], Section 2.
\bigskip

For technical reasons we start with a situation dual to the above one. So let
$F$, $G$ be finite free $R$-modules of rank $m,n$ and $\psi: G \to F$ an
$R$-homomorphism. Set $\hat G=G\otimes \SS(F)$ where $\SS(F)$ denotes the
symmetric algebra of $F$. Then we may consider $\psi$ an $\SS(F)$-linear form
on $\hat G$ and can define the Koszul antiderivation
$$
\partial_\psi:\bigwedge\hat G\longrightarrow \bigwedge\hat G
$$
with respect to $\psi$ in the usual way, i.e.
$$
\partial_\psi(x_1\wedge\ldots\wedge x_i)=\sum_{j=1}^i(-1)^{j+1}\psi(x_j)x_1
\wedge \ldots{\widehat x}_j\ldots\wedge x_i $$ for $x_1\ldots x_i\in \hat G$.
We use the term Koszul complex also for the complex $$ 0\to
R\stackrel{\phi}{\to}G\stackrel{\phi(1)\wedge}{\to}\bigwedge^2
G\stackrel{\phi(1)\wedge}{\to}\bigwedge^3 G\stackrel{\phi(1)\wedge}{\to}\dots
$$ associated with a linear map $\phi:R\to G$. Suppose that $\psi\phi=0$ and
let $$ d_\phi:\bigwedge\hat G\longrightarrow \bigwedge\hat G $$ be the
differential of the Koszul complex associated with $\phi\otimes \SS(F)$, i.~e.
$$ d_\phi(x)=(\phi(1)\otimes 1)\wedge x $$ for $x\in \bigwedge\hat G$. Since
$\psi\phi=0$, obviously $\partial_\psi d_\phi+d_\phi\partial_\psi =0$. Thus we
obtain the \emph{Koszul bicomplex} $\KK$ $$
\begin{CD}
&&0&&0&&0&&0&&0&&\\
&& @VVV @VVV @VVV @VVV @VVV \\
0@>>> R@>\phi>> G@>>> \cdots  \bigwedge^{p-1} G@>d_\phi>>\bigwedge^p G
   @>>> \bigwedge^{p+1} G && \cdots\\
&& @VVV @V\psi VV @VVV @V\partial_\psi VV @VVV \\
&&0@>>> F@>>> \cdots M^{p-1,1}@>>> M^{p,1}@>>> M^{p+1,1}&& \cdots\\
&&&& @VVV @VVV @VVV @VVV \\
&&&&0@>>> \cdots  M^{p-1,2}@>>> M^{p,2}@>>> M^{p+1,2}&& \cdots\\
&&&&&& @VVV @VVV @VVV\\
&&&&&& \vdots && \vdots && \vdots\\
&&&&&& 0@>>> M^{p,p}@>>> M^{p+1,p}&& \cdots\\
&&&&&&&& @VVV @VVV \\
&&&&&&&& 0@>>> M^{p+1,p+1} &&\cdots \\
&&&&&&&&&& @VVV\\
&&&&&&&&&& 0 &&\cdots
\end{CD}
$$
where
$$
M^{p,q}=\bigwedge^{p-q}G\otimes\SS^q(F)
$$
for all integers $p$, $q$, and $\SS^q$ means $q$th symmetric power. The row
homology of $\KK$ at $M^{p,q}$ is denoted by $H_\phi^{p,q}$, the column
homology by $H_\psi^{p,q}$. Thus $H_\phi^{p,0}$ is the $p$th homology module
$H^p$ of the Koszul complex associated to $\phi$. Set $N^p= \Ker
\partial_\phi^{p,0}$. The canonical injections $N^p\to \bigwedge^p G$ yield a
complex homomorphism
$$
\begin{CD}
0@>>> R@>\bar\phi>>N^1 @>>>\cdots && N^p @>d_{\bar\phi}>> N^{p+1}&&\cdots \\
&& \parallel && @VVV && @VVV @VVV \\
0@>>> R@>\phi>> G@>>>\cdots && \bigwedge^p G@>d_\phi >> \bigwedge^{p+1}G&&\cdots
\end{CD}
$$ where the maps $\bar\phi$, $d_{\bar\phi}$ are induced by $\phi$, $d_\phi$.
The homology of the first row at $N^p$ is denoted by $\bar{H}^p$. We are now
ready to state the key proposition. Here as in the following $^*$ means
$\Hom_R(\;\; ,R)$. Moreover, $\I_\Psi$ denotes the ideal of $m$-minors of (a
matrix representing) $\psi$.

\begin{proposition} Set $g=\grade\I_M$, $C=\Cok\psi$, and let $h=\grade\Im
\phi^*$.  Assume that $r=n-m\ge 1$ and $g=r+1$. Then
\begin{enumerate}
\item[\rm{(a)}] $\Im\phi^*\subset \I_M$ and, in particular, $h\le g$;
\item[\rm{(b)}]  $\bar H^i=0$ for $0\le i\le\min(2,h-1)$;
\item[\rm{(c)}]
$$ \bar{H}^i= \begin{cases} \SS^{\frac{i-1}2}(C) &\quad\textrm{ if }\quad 3\le
i< h,\; i\not\equiv 0\ (2)\\ 0 &\quad\text{if}\quad 3\le i< h,\; i\equiv 0\
(2);
\end{cases}
$$
\item[\rm{(d)}] moreover, for $h\ge 3$ there is an exact sequence
\begin{alignat*}{2}
&0\to \SS^{\frac{h-1}2}(C)\to \bar{H}^h \to H^h\quad&&
\text{if } h\not\equiv 0\ (2),\\
&0\to \bar{H}^h\to  H^h&&\text{if } h\equiv 0\ (2).
\end{alignat*}
\end{enumerate}
\end{proposition}

\begin{proof} Let $M=\Cok\psi^*$. We choose a basis $e_1,\dots,e_m$ of $F^*$ and
define the linear map $\Psi:G^*\to \bigwedge^{m+1}G^*$ by
$\Psi(x)=\psi^*(e_1)\wedge\ldots\wedge\psi^*(e_m)\wedge x$. Then one obtains a
complex
$$
0\to F^*\stackrel{\psi^*}{\to} G^*\stackrel{\Psi}{\to} \bigwedge^{m+1}G^*
$$
whose dual is the head of the Buchsbaum--Rim complex resolving $C=\Cok \psi$.
It follows that $M^*=\Im \Psi^*$, and obviously $\Im\Psi^*\subset \I_M G$.
Since $\phi\in M^*$, one has $\Im\phi\subset \I_M$. This shows (a).

We quote some well-known facts about the homology of $\KK$. Let $0\le p\le g$.
Then $H_\psi^{p,q}=0$ for $q\ne 0,p$ and $H_\psi^{p,p}=\SS^p(C)$. (See [BV1],
Proposition 2.1 for the general statement.) Furthermore $H_\phi^{p,q}=0$ for
$p-q<h$ by the grade sensitivity of the Koszul complex for $\phi$.

Claim (b) on $\bar H^i$ for $0\le i\le \min(2,h-1)$ is easily proved from the
long exact (co)homology sequence.

For (c) we modify the complex $\KK$ to the complex $\tilde\KK$ by setting (i)
$M^{p,p+1}=\SS^p(C)$ and (ii) $M^{p,-1}=N^p$. The maps to be added are the
natural surjection $M^{p,p}\to \SS^p(C)$, the zero map $\SS^p(C)\to
M^{p+1,p+1}$, and those induced by the canonical injections $N^p\to\bigwedge^p
G$. The truncation of $\tilde \KK$ to the ``rectangle'' $0\le p\le h$, $0\le
q\le g$ has exact columns. Moreover, row homology for indices $<h$ can only
occur at $M^{p,-1}$, namely $\bar H^p$, $0\le p\le h$, and at $M^{p,p+1}$,
namely $\SS^p(C)$.

For an inductive argument we let $\RR_q$, $q\ge -1$, be the $q$th row of
$\tilde \KK$ and $\BB_{q+1}$ be the image complex of $\RR_q$ in $\RR_{q+1}$.
Then we have a series of exact sequences
$$
0\to \BB_q\to \RR_q \to \BB_{q+1}\to 0.
$$
Thus we can use the long exact (co)homology sequence for each $q$. With
$E^{p,q}=H^p(\BB_q)$ one therefore obtains the ``southwest'' isomorphisms
$$
\bar H^i=E^{i,0}\iso E^{i-1,1}\iso\dots\iso E^{\frac i2,\frac i2}
$$
if $i$ is even, and
$$
\bar H^i=E^{i,0}\iso E^{i-1,1}\iso\dots\iso E^{\frac {i+1}2,\frac {i-1}2}
$$
if $i$ is odd, $0\le i<h$. In fact, there is an exact sequence
$$
H^{i-(j+1),j} \to E^{i-(j+1),j+1}\to E^{i-j,j} \to H^{i-j,j}
$$
and the extreme terms in this sequence are $0$ for all $j$ under consideration.

Let now $i$ be even, and $j=i/2$. Since the map $M^{j,j}\to M^{j+1,j}$ is
injective, the same holds true for its restriction $\BB^j_j\to \BB^{j+1}_j$ in
$\BB_j$, and so $\bar H^p=E^{j,j}=0$.

In the case in which $i$ is odd we can go one further step southwest, and
obtain the isomorphism $\SS^j(C)=E^{j,j+1}\iso E^{j+1,j}$ for $j=(i-1)/2$.

Since $H^h\neq 0$, we only have an exact sequence
$$
0\to E^{h-1,h}\to \bar H^h\to H^h,
$$
but the arguments above can be applied to $E^{h-1,1}$; it is zero if $h$ is
even, and isomorphic to $\SS^{(h-1)/2}(C)$ if $h$ is odd.--
\end{proof}

\begin{remark} \textrm{Proposition 1} shows that roughly the first half of the
symmetric powers $\SS^i(C)$, $i\le h$, can be interpreted as homologies of a
Koszul complex. It is also possible to interpret the other half in a similar
way. To this end we consider the column complexes $\CC_0,\dots,\CC_h$ of $\KK$
({\it not} of $\tilde\KK$) and set $\CC_{h+1}=\Cok(\CC_{h-1}\to\CC_h)$. Then we
obtain an exact sequence
$$
0\to\CC_0\to\dots\to\CC_h\to\CC_{h+1}\to 0
$$
of complexes, and the only nonzero (co)homology can occur along $\CC_{h+1}$ and
at $H^p(\CC_p)\iso S_p(C)$ for $p>0$. If one applies arguments similar to those
in the proof of Proposition 1 (now proceeding in northwestern direction), then
one obtains
$$
H^i(\CC_{h+1})\iso\begin{cases}
 \SS^{\frac{h+i}2}(C)& \quad\text{if } h+i\equiv 0\ (2),\\
 0                   & \quad\text{if } h+i\not\equiv 0\ (2).
\end{cases}
$$
for $0\le i\le h$. Note that the module $\CC_{h+1}^i$ is resolved by $\RR_i$,
and $\RR_i$ is just a truncated and shifted version of $\RR_0\tensor \SS^i(F)$.
The truncations of $\RR_0$ resolve the exterior powers $\bigwedge^j M_\phi$
where $M_\phi=\Cok(\phi)$. Thus $\CC_{h+1}^i=\bigwedge^{h-i} M_\phi\tensor
\SS^i(F)$.--
\end{remark}

As in the proof of Proposition 1 let $M=\Cok\psi^*$. Then it is easy to see
that $N^p=(\bigwedge^p M)^*$ for all $p$. In fact, $\psi^*$ induces a
presentation
$$
\bigwedge^{p-1}G^*\otimes F^*\to \bigwedge^p G^*\to \bigwedge^p M\to 0
$$
for all $p$, and the exact sequence $0\to N^p\to\bigwedge^p G\to
\bigwedge^{p-1}G\otimes F$ is obtained by dualizing. It
follows that $N^r$ is free of rank $1$ (provided $\grade \I_M\ge 2$), and
$N^p=0$ for $p>r$.

\begin{corollary} Let $\psi:G\to F$ be as above, and assume that
$g=\grade\I_M=r+1$. Then the following conditions are equivalent.
\begin{enumerate}
\item[\textrm{(1)}] There is a homomorphism $\phi:R\to G$ such that
$\psi\phi=0$ and the ideal $\Im\phi^*$ has grade $r+1$;
\item[\textrm{(2)}] {\rm (i)} $r=1$ or {\rm (ii)} $m=1$ and $r$ is odd.
\end{enumerate}
\end{corollary}

\begin{proof} The implication $(2)\Rightarrow (1)$ is an easy exercise. (See
also the considerations at the end of this note.)

For the converse observe that $N^p=0$ for $p>r$ and that $N^r$ is free of rank
$1$. So $\bar{H}^r$ must be cyclic. If $r$ were even, then $\bar{H}^r=0$ by
Proposition 1, and $\Im\phi^*=R$. Thus $r$ must be odd. In this case
$\bar{H}^r=\SS^{(r-1)/2}(C)$ by Proposition 1 where $C=\Cok\psi^*$. So if $r\ge
3$, then $C$ must be cyclic, which in turn means $m=1$.--
\end{proof}

We now return to our original purpose. Therefore let $\bar\chi$ be a linear
form on $M=\Cok\psi^*$. The induced linear form on $G^*$ is denoted by $\chi$;
note that  $\chi\psi^*=0$. Set $\phi=\chi^*$, and, as above, $r=n-m$. We want
to connect the truncated Koszul complex
\begin{equation}
0 \to \bigwedge^r M \to \cdots \bigwedge^i M\stackrel{\partial_{\bar\chi}}{\to}
 \bigwedge^{i-1} M\cdots \to M \to R\to 0\tag{1}
\end{equation}
with the complex
\begin{equation}
0\to R\stackrel{\bar\phi}{\to}N^1 \to\cdots N^p\stackrel{d_{\bar\phi}}{\to}
N^{p+1}\cdots\tag{2}
\end{equation}
considered above. We have already observed that $N^p=(\bigwedge^p M)^*$.

\begin{lemma} With notation from above, let $g=\grade\I_M\ge r+1$. Then there are maps $\mu_p:\bigwedge^p
M\to N^{r-p}$, $p=0,\ldots, r$, such that $d_{\bar\phi}\mu_p=\pm
\mu_{p-1}\partial_{\bar\chi}$ and which are isomorphisms for $p=0,\ldots, r-1$
and injective for $p=r$. If, in particular, $g=r+1\ge 2$, then we obtain the
following diagram of maps, the columns of which are exact and with commutative
or anticommutative rectangles: $$
\begin{CD}
&& && 0 &&  0\\
&& && @VVV  @VVV\\
0@>>>\bigwedge^r M@>>>\bigwedge^{r-1} M@>>> \bigwedge^{r-2} M\\
&& @V\mu_r VV  @V\mu_{r-1} VV @VVV \\
0 @>>> R@>\bar\phi>>N^1 @>>> N^2&&\\
&& @VVV  @VVV \\
0@>>> R/\I_M@>>>  0\\
&& @VVV\\
&& 0
\end{CD}
$$
\end{lemma}

\begin{proof} As in [BV2], p.~17, we choose isomorphisms $\gamma:\bigwedge^m
F^*\to R$,  $\delta:\bigwedge^n G^*\to R$, and define maps
$$
\nu_p:\bigwedge^p G^* \to (\bigwedge^{r-p} G^*)^*,
$$
$p=0,\ldots, r$, by
$$
\nu_p(x)(y)=\delta(x\wedge y\wedge\bigwedge^m\psi^*(z)),
$$
$x\in \bigwedge^p G^*,\ y\in\bigwedge^{r-p} G^*,\ z=\gamma^{-1}(1)$. Via the
natural isomorphism $(\bigwedge^{r-p} G^*)^*\cong \bigwedge^{r-p} G$ we regard
$\nu_p$ as a map $\bigwedge^p G^* \to \bigwedge^{r-p} G$. One has
$\Im\nu_p\subset N^{r-p}$, and it is an easy exercise to show that the diagram
$$
\begin{CD}
\bigwedge^p G^*@>\partial\chi>> \bigwedge^{p-1} G^* \\
@V\nu_pVV @V\nu_{p-1}VV\\
\bigwedge^{r-p} G@>d_{\chi^*}>> \bigwedge^{r-p+1} G
\end{CD}
$$
is commutative or anticommutative (see for example [HM], proof of Theorem 3.1).
Consequently the same is true for
$$
\begin{CD}
\bigwedge^p M @>\partial\bar\chi>>\bigwedge^{p-1} M\\
@V\rho_pVV @V\rho_{p-1}VV\\
\Im\nu_p@>{d_{\chi^*}|\Im\nu_p}>>\Im\nu_{p-1}
\end{CD}
$$
where $\rho_p$ and $\rho_{p-1}$ are induced by $\nu_p$ and $\nu_{p-1}$. Now let
$\mu_p$ be the composition of $\rho_p$ and the canonical injection $\Im\nu_p\to
N^{r-p}$. Then the equation asserted in the proposition obviously holds. In
case $p<\grade\I_M-1$,
$\Im\nu_p=N^{r-p}$. This proves the remaining statements.--
\end{proof}

If we look at the homology of the truncated Koszul complex (1) associated to
$\bar\chi$ instead of the homology of (2), then a somewhat smoother assertion
than Proposition 1 is possible.

\begin{theorem} Let $M$ be module with a finite free resolution of length
$1$, $M=\Cok\psi^*$ where $\psi:G\to F$ is as above, and assume that
$g=\grade\I_M =r+1$. Let $\bar\chi$ be a linear form on $M$. Then
$\Im\bar\chi\subset \I_M$, and for the homology $\bar H_p$ at $\bigwedge M^p$
of the truncated Koszul complex (1) associated to $\bar\chi$ the following
holds:
$$
\bar H_{r-i}=
\begin{cases}
 \SS^{\frac{i-1}2}(C) &\quad\text{if}\ \ 0\le i< h,\; i\not\equiv 0\ (2)\\
0 &\quad\text{if}\ \ 0\le i< h,\; i\equiv 0\ (2),
\end{cases}
$$
where $\SS^0(C)=R/\I_M$.

Furthermore there is a $\bar\chi$ with $\grade\Im\bar\chi=g$ if and only if
{\rm (i)} $r=1$ or {\rm (ii)} $m=1$ and $r$ is odd. In this case we have
necessarily $\Im\chi=\I_M$.
\end{theorem}

\begin{proof} Consider the diagram in Lemma 4. Since $\mu_r$ is injective,
$\bar H_r$ must be zero if $h>0$. Next let $h>1$. Then we obtain $R/\I_M\cong
\bar H_{r-1}$, since $\bar{H}^0=\bar{H}^1=0$. Finally if $h>2$, then in
addition $\bar{H}^2=0$, so $\bar H_{r-2}=0$ as stated. The remaining assertions
concerning the homology of (1) are contained in Proposition 1.

Instead of $\bar\chi$ we can consider the induced linear form $\chi$ on $G^*$.
Corollary 3 yields the statement about the existence of such a linear form
$\chi$ satisfying $\grade\Im\chi=g$. Assume that such a $\chi$ exists. If
$m=1$, then $\Im\chi=\I_M$ by Proposition 2 in [BV3]. If $r=1$, then, under our
assumptions, $M=\Cok\psi^*$ is an ideal in $R$ which must be isomorphic to
$\Im\chi$. Using the Hilbert-Burch Theorem, we have $\Im\chi=a\I_M$ with an
element $a\in R$. Since $\grade\Im\chi=2$, $a$ must be a unit.--
\end{proof}

\begin{remark}
\begin{enumerate}
\item[(a)] It is a noteworthy fact that the homology $\bar
H_{r-i}$ of the truncated Koszul complex (1) associated to $\bar\chi$ is
independent of $\chi$ for $i<h$.

\item[(b)] The Koszul complex associated to a linear form $\chi$ on a free
module of rank $r$ is grade sensitive: its homology vanishes for $j>r-\grade
\Im\chi$, and does not vanish at $r-\grade\Im\chi$.

In a sense, this is also true for the linear form $\bar\chi$ considered in
Theorem 5: of course, ``vanishes'' must be replaced by ``vanishes if $i$ is
even and is isomorphic to $\SS^{(i-1)/2}(C)$ if $i$ is odd''. Then Theorem 5
covers all $i=0,\dots,h-1$, but the analogy also persists if $i=h<g$. In fact,
let $\pp$ be a prime ideal of grade $h$ such that $\Im \chi\subset\pp$. The
module $M_\pp$ is free and $\SS^j(C_\pp)=\SS^j(C)_\pp=0$ for all $j$. Therefore
 one can apply the grade sensitivity of the Koszul complex for a free module,
and $\bar H_h$ can be neither $0$ nor isomorphic to a non-vanishing symmetric
power of $C$: otherwise its localization would vanish.

\item[(c)] That we have truncated the Koszul complex of $\bar\chi$ is
inessential. In fact, $\bigwedge^r M$ is torsionfree of rank $1$, and
$\bigwedge^{r+1} M$ has rank $0$. Hence $\Hom_R(\bigwedge^{r+1} M, \bigwedge^r
M)=0$, and the homology of the full and of the truncated Koszul complexes at
$\bigwedge^r M$ coincide.
\end{enumerate}
\end{remark}

Let $\psi:G\to F$ be as above, $r=n-m$, and $g=\grade\I_M=r+1$. In case
$r>1$, the existence of a linear form $\chi$ on $G^*$ with $\grade\Im\chi = g$
can be described equivalently and independently of the last theorem.

\begin{proposition} With the assumptions on $\psi$ and $\chi$ from Theorem
5, assume in addition that $r>1$. Then $\grade\Im\chi= g$ is possible if and
only if there exists a submodule $U$ of $M=\Cok\psi^*$ with the following
properties:
\begin{enumerate}
\item[\textrm{(1)}] $\rank U=r-1$;
\item[\textrm{(2)}] $U$ is reflexive, orientable, and $U_{\frak p}$ is a free
direct summand of $M_{\frak p}$ for all primes $\frak p$ of $R$ such that
$\grade\frak p\le r$.
\end{enumerate}
\end{proposition}

\begin{proof}
Let $\bar\chi$ be a linear form on $M$ such that
$\grade\Im\bar\chi= r+1$. Set $U=\Ker\bar\chi$. Then (1) and the last property
in (2) are obvious. Since $\Im\bar\chi$ is torsionfree and $M$ is reflexive,
$U$ must be reflexive. Because $\Im\bar\chi$ has grade $\ge 2$, it is
orientable. $M$ being orientable, the orientability of $U$ follows from
Proposition (2.8) in [B].

Conversely let $U$ be a submodule of $M$ which satisfies (1) and (2). Then
$I=M/U$ is torsionfree of rank 1 and therefore an ideal in $R$ which is
orientable since $U$ and $M$ are orientable. Consequently $\grade I\ge 2$. So
for a prime $\frak p$ in $R$ which contains $I$, the localization $IR_{\frak
p}$ cannot be free. On the other hand $I=R$ is impossible since $g=r+1$. In
view of the last condition in (2), $I$ must have grade $r+1$.--
\end{proof}

>From Theorem 5 we know that the hypothesis of Proposition 7 can only be
satisfied with $m=1$ and $r$ odd. The submodule $U$ of $M$ in Proposition 7 has
projective dimension $r-1$. In fact, the ideal $\I_M=\Im \chi$ is generated by
$r+1$ elements and has grade $r+1$. Therefore it is perfect, i.~e.\ $\projdim
R/I=r+1$. This implies $\projdim U=r-1$ via the exact sequence $0\to U\to M\to
I\to 0$. Dualizing this exact sequence, we obtain an exact sequence $0\to R \to
M^* \to U^*\to 0$. Since $M^*\cong \bigwedge^{r-1}M$ has projective dimension
$r-1$, it follows that $U^*$ has projective dimension $r-1$.

The dualization argument just used amounts to interchanging the roles of
$\psi^*$ and $\chi$, so that $U^*$ plays the same role for $\chi^*$ and $\psi$
as $U$ for $\psi^*$ and $\chi$. If we choose $\chi$ in a special way, then we
can even achieve that $U$ and $U^*$ are isomorphic in a skewsymmetric way: for
the isomorphism $\sigma:U\to U^*$ one has $\sigma^*=-\sigma$ upon the
identification of $U$ and $U^{**}$ via the natural isomorphism.

As we mentioned at the beginning of the proof of Corollary 3, it is easy to see
that there is a linear form $\chi$ on $G^*$ such that $\psi^*(1)\in \Ker\chi$
and $\grade\Im\chi=r+1=n$: fix a basis $z_1,\ldots,z_n$ of $G^*$ and let
$\psi^*(1)=\sum_{i=1}^n x_iz_i$; the map $\chi: \sum_{i=1}^n a_iz_i \mapsto
\sum_{i=1}^n(-1)^ia_ix_{n+1-i}$ is an appropriate linear form. Let $\bar\chi$
be the induced form on $M=\Cok\psi^*$. The submodule $U=\Ker\bar\chi$ has the
properties (1) and (2) of Proposition 7 (see the first part of the proof).
Consider the commutative diagram
$$
\begin{CD}
0@>>> \Ker\chi @>>> G^* @>\chi>> R\\
&& @V\rho_1 VV @V\rho VV \parallel\\
0@>>> N @>>> G @>\psi >> R,
\end{CD}
$$
with exact rows. The isomorphism $\rho$ is defined by
$\rho(z_i)=(-1)^iz^*_{n+1-i}$ where $z_1^*,\ldots,z_n^*$ is the basis of $G$
dual to $z_1,\ldots,z_n$, and $\rho_1$ is  induced by $\rho$. Since
$\rho_1(\psi^*(1))=-\bar\chi$, we obtain a second commutative diagram
$$
\begin{CD}
0@>>> R^* @>(\psi^*)_1>> \Ker\chi @>>> U @>>> 0\\
&& @VVV @V\rho_1 VV @V\bar\rho_1VV\\
0@>>> R^* @>\bar\chi^*>> N @>>> \phantom{*}U^* @>>> 0
\end{CD}
$$
with exact rows. $(\psi^*)_1$ is induced by $\psi^*$ and the first vertical
arrow maps $1$ to $-1$. As $\rho_1$ is an isomorphism, the same is true for
$\bar\rho_1$, so $U$ turns out to be self-dual. Moreover $\bar\rho_1$ is
skew-symmetric, i.e. $(\bar\rho_1)^*=-\bar\rho_1$ since the same is true for
$\rho$.

Suppose that $R=K[X_1,\dots,X_n]$ is the polynomial ring in $n$ indeterminates
over a field $K$. If we define $\psi:R^n\to R$ by $\psi(e_i)=X_i$, then the
module $U$ is associated with a rank $n-2$ vector bundle on $\PP^{n-1}(K)$.
Such bundles exist however also for odd $n$; see [V]. The module $V$
assoc\'{\i}ated with the the construction in [V] is self-dual only for $n=4$.

\end{document}